\documentclass{elsart}
\usepackage{amsmath}
\usepackage{amssymb}
\usepackage{graphicx}

\newcommand{\R}{\mathbb{R}}
\newcommand{\Cdot}{\boldsymbol{\cdot}}

\begin{document}

\begin{frontmatter}

\title{Distribution functions of linear combinations of lattice polynomials from the uniform distribution}

\author{Jean-Luc Marichal}
\ead{jean-luc.marichal[at]uni.lu}

\address{
Institute of Mathematics, University of Luxembourg\\
162A, avenue de la Fa\"{\i}encerie, L-1511 Luxembourg, Luxembourg }

\author{Ivan Kojadinovic\thanksref{Ivan}}
\ead{ivan[at]stat.auckland.ac.nz}
\thanks[Ivan]{Formerly at LINA CNRS FRE 2729, \'Ecole Polytechnique de l'Universit\'e de Nantes, France.}

\address{
Department of Statistics, The University of Auckland, Private Bag 92019 \\
Auckland 1142, New Zealand}

\date{Submitted 18 Jul 06, Accepted 15 Sep 07}

\begin{abstract}
We give the distribution functions, the expected values, and the moments of linear combinations of lattice polynomials from the uniform
distribution. Linear combinations of lattice polynomials, which include weighted sums, linear combinations of order statistics, and lattice
polynomials, are actually those continuous functions that reduce to linear functions on each simplex of the standard triangulation of the unit
cube. They are mainly used in aggregation theory, combinatorial optimization, and game theory, where they are known as discrete Choquet
integrals and Lov\'asz extensions.
\end{abstract}

\end{frontmatter}

%---------------------------------------------------------------------------------------------- Section 1
\section{Introduction}

Let $h:[0,1]^n\to \R$ be an aggregation function and let $\mathbf{X}$ be a random vector uniformly distributed on $[0,1]^n$. An interesting but
generally difficult problem is to provide explicit expressions for the distribution function and the moments of the aggregated random variable
$Y=h(\mathbf{X})$.

This problem has been completely solved for certain aggregation functions (see for instance\ \cite[\S 7.2]{Mcc04}), especially piecewise linear
functions such as weighted sums \cite{BarSmi79} (see also \cite{MarMos}), linear combinations of order statistics \cite{AgaDalSin02,Mat85,Wei71}
(see also \cite[\S 6.5]{DavNag03} for an overview), and lattice polynomials \cite{Mar06}, which are max-min combinations of the variables.

In this note we solve the case of linear combinations of lattice polynomials, which include the three above-mentioned cases. Actually, linear
combinations of lattice polynomials are exactly those continuous functions that reduce to linear functions on each simplex of the standard
triangulation of $[0,1]^n$. In particular, these functions are completely determined by their values at the $2^n$ vertices of $[0,1]^n$.

The concept of linear combination of lattice polynomials is known in combinatorial optimization and game theory as the {\em Lov\'asz
extension}\/ \cite{AlgBilFerJim04,GraMarRou00,Lov83,Sin84} of a pseudo-Boolean function (recall that a pseudo-Boolean function is a real-valued
function of 0-1 variables). When it is nondecreasing in each variable, it is known in the area of nonlinear aggregation and integration as the
discrete {\em Choquet integral}\/ \cite{Den94,GraMurSug00,Mar02b}, which is an extension of the discrete Lebesgue integral (weighted mean) to
non-additive measures. The equivalence between the Lov\'asz extension and the Choquet integral is discussed in \cite{Mar02b}.

This note is set out as follows. In Section 2 we elaborate on the definition of linear combinations of lattice polynomials and we show how to
concisely represent them. In Section 3 we provide formulas for the distribution function and the moments of any linear combination of lattice
polynomials from the uniform distribution. Finally, in Section 4 we provide an application of our results to aggregation theory.

Throughout we will use the notation $[n]:=\{1,\ldots,n\}$. Also, for any subset $A\subseteq [n]$, $\mathbf{1}_A$ will denote the characteristic
vector of $A$ in $\{0,1\}^n$. Finally, for any function $h:[0,1]^n\to\R$, we define the set function $v_h:2^{[n]}\to\R$ as
$v_h(A):=h(\mathbf{1}_A)$ for all $A\subseteq [n]$.

%---------------------------------------------------------------------------------------------- Section 2
\section{Linear combinations of lattice polynomials}

In the present section we recall the definition of lattice polynomials and we show how an arbitrary combination of lattice polynomials can be
represented.

Basically an $n$-place lattice polynomial $p:[0,1]^n\to [0,1]$ is a function defined from any well-formed expression involving $n$ real
variables $x_1,\ldots,x_n$ linked by the lattice operations $\wedge=\min$ and $\vee=\max$ in an arbitrary combination of parentheses (see e.g.\
Birkhoff \cite[\S II.2]{Bir67}). For instance,
$$
p(x_1,x_2,x_3)=(x_1\wedge x_2)\vee x_3
$$
is a 3-place lattice polynomial.

Consider the standard triangulation of $[0,1]^n$ into the canonical simplices
\begin{equation}\label{eq:Ssigma}
S_{\sigma}:=\{x\in [0,1]^n\mid x_{\sigma(1)}\geqslant\cdots\geqslant x_{\sigma(n)}\}\qquad (\sigma\in\frak{S}_n),
\end{equation}
where $\frak{S}_n$ is the set of all permutations on $[n]$. Clearly, any linear combination of $n$-place lattice polynomials
$$
h(\mathbf{x})=\sum_{i=1}^m c_i\, p_i(\mathbf{x})%\qquad (\mathbf{x}\in [0,1]^n)
$$
is a continuous function whose restriction to any canonical simplex is a linear function. According to Singer \cite[\S 2]{Sin84}, $h$ is then
the {\em Lov\'asz extension}\/ of the pseudo-Boolean function $h|_{\{0,1\}^n}$, that is, the continuous function defined on each canonical
simplex $S_{\sigma}$ as the unique linear function that coincides with $h|_{\{0,1\}^n}$ at the $n+1$ vertices
$$
\varepsilon_i^{\sigma}:=\mathbf{1}_{\{\sigma(1),\ldots,\sigma(i)\}} \qquad (i=0,\ldots,n)
$$
of $S_{\sigma}$. It can be written as \cite[\S 2]{Sin84}
\begin{equation}\label{eq:LovExt}
h(\mathbf{x})=\sum_{i=1}^n \big(h_i^{\sigma}-h_{i-1}^{\sigma}\big)\, x_{\sigma(i)} \qquad (\mathbf{x}\in S_{\sigma}),
\end{equation}
where $h_i^{\sigma}:=h\big(\varepsilon_i^{\sigma}\big)=v_h\big(\{\sigma(1),\ldots,\sigma(i)\}\big)$ for all $i=0,\ldots,n$. In particular,
$h_0^{\sigma}=0$.

Conversely any continuous function $h:[0,1]^n\to\R$ that reduces to a linear function on each canonical simplex is a linear combination of
lattice polynomials:
\begin{equation}\label{eq:MoEx}
h(\mathbf{x})=\sum_{A\subseteq [n]}m_h(A)\,\bigwedge_{i\in A}x_i \qquad (\mathbf{x}\in [0,1]^n),
\end{equation}
where $m_h:2^{[n]}\to\R$ is the {\em M\"obius transform}\/ of $v_h$, defined as
$$
m_h(A):=\sum_{B\subseteq A}(-1)^{|A|-|B|}v_h(B).%\qquad (A\subseteq [n]).
$$
Indeed, expression (\ref{eq:MoEx}) reduces to a linear function on each canonical simplex and agrees with $h(\mathbf{1}_B)$ at $\mathbf{1}_B$
for each $B\subseteq [n]$.

Eq.~(\ref{eq:LovExt}) thus provides a concise expression for linear combinations of lattice polynomials. We will use it in the next section to
calculate their distribution functions and their moments.

%\begin{exmp}
%Consider the function
%$$
%h(\mathbf{x})=2[x_1\wedge(x_2\vee x_3)]-3(x_1\wedge x_3)+x_1,
%$$
%for which we have $v_h(\{1\})=1$, $v_h(\{1,2\})=3$, and $v_h(A)=0$ for $A\neq \{1\}$ and $A\neq \{1,2\}$. According to (\ref{eq:MoEx}), it can
%also be rewritten as
%$$
%h(\mathbf{x})=x_1+2(x_1\wedge x_2)-(x_1\wedge x_3)-2(x_1\wedge x_2\wedge x_3).
%$$
%\end{exmp}

\begin{rem}\label{rem:1}
As we have already mentioned, the class of linear combinations of lattice polynomials covers three interesting particular cases, namely: lattice
polynomials, linear combinations of order statistics, and weighted sums. These are characterized as follows. Let $h:[0,1]^n\to\R$ be a linear
combination of lattice polynomials.
\begin{enumerate}
\item The function $h$ reduces to a lattice polynomial if and only if the set function $v_h$ is monotone, $\{0,1\}$-valued, and such that
$v_h([n])=1$.

\item As the order statistics are exactly the symmetric lattice polynomials (see \cite{Mar02c}), the function $h$ reduces to a linear
combination of order statistics if and only if the set function $v_h$ is cardinality-based, that is, such that $v_h(A)=v_h(A')$ whenever
$|A|=|A'|$.

\item The function $h$ reduces to a weighted sum if and only if the set function $v_h$ is additive, that is, $v_h(A)=\sum_{i\in A}v_h(\{i\})$.
\end{enumerate}
\end{rem}

%---------------------------------------------------------------------------------------------- Section 3
\section{Distribution functions and moments}

Before yielding the main results, let us recall some basic material related to divided differences. See for instance \cite{Dav75,DeVLor93,Pow81}
for further details.

Consider the {\em plus}\/ (resp.\ {\em minus}) {\em truncated power function}\/ $x^n_+$ (resp.\ $x^n_-$), defined to be $x^n$ if $x>0$ (resp.\
$x<0$) and zero otherwise. Let $\mathcal{A}^{(n)}$ be the set of $n-1$ times differentiable one-place functions $g$ such that $g^{(n-1)}$ is
absolutely continuous. The $n$th {\em divided difference}\/ of a function $g\in \mathcal{A}^{(n)}$ is the symmetric function of $n+1$ arguments
defined inductively by $\Delta[g:a_0]:=g(a_0)$ and
$$
\Delta[g:a_0,\ldots,a_n]:=
\begin{cases}
\displaystyle{\frac{\Delta[g:a_1,\ldots,a_n]-\Delta[g:a_0,\ldots,a_{n-1}]}{a_n-a_0}}, & \mbox{if $a_0\neq a_n$},\\
\displaystyle{\frac{\partial}{\partial a_0}\,\Delta[g:a_0,\ldots,a_{n-1}]}, & \mbox{if $a_0= a_n$.}
\end{cases}
$$

The {\em Peano representation}\/ of the divided differences, which can be obtained by a Taylor expansion of $g$, is given by %\cite{CurSch66}
\begin{equation}\label{eq:Peano}
\Delta[g:a_0,\ldots,a_n]=\frac 1{n!}\,\int_{\R} g^{(n)}(t)\, M(t\mid a_0,\ldots,a_n)\,\mathrm{d}t,
\end{equation}
where $M(t\mid a_0,\ldots,a_n)$ is the {\em B-spline} of order $n$, with knots $\{a_0,\ldots,a_n\}$, defined as
\begin{equation}\label{eq:BSpline}
M(t\mid a_0,\ldots,a_n):=n\,\Delta[(\Cdot-t)^{n-1}_+:a_0,\ldots,a_n].
\end{equation}
%with support $\mathrm{conv}\{a_0,\ldots,a_n\}$.

We also recall the {\em Hermite-Genocchi formula}: For any function $g\in \mathcal{A}^{(n)}$, we have
\begin{equation}\label{eq:HermGeno}
\Delta[g:a_0,\ldots,a_n]=
%\int_{[a_0,\ldots,a_n]}g^{(n)}:=
\int_{S_{id}}g^{(n)}\Big[a_0+\sum_{i=1}^n (a_i-a_{i-1})x_i\Big]\,\mathrm{d}\mathbf{x},
\end{equation}
where $S_{id}$ is the simplex defined in (\ref{eq:Ssigma}) when $\sigma$ is the identity permutation.

For distinct arguments $a_0,\ldots,a_n$, we also have the following formula, which can be verified by induction,
\begin{equation}\label{eq:DvDDistArg}
\Delta[g:a_0,\ldots,a_n]=\sum_{i=0}^n\frac{g(a_i)}{\prod_{j\neq i}(a_i-a_j)}.
\end{equation}

Now, consider a random vector $\mathbf{X}$ uniformly distributed on $[0,1]^n$ and set $Y_h:=h(\mathbf{X})$, where the function $h:[0,1]^n\to\R$
is a linear combination of lattice polynomials as given in formula (\ref{eq:LovExt}). We then have the following result.

\begin{thm}\label{thm:ExpgY}
For any function $g\in \mathcal{A}^{(n)}$, we have
\begin{equation}\label{eq:ExpgY}
\mathbf{E}[g^{(n)}(Y_h)]=\sum_{\sigma\in\frak{S}_n}\Delta[g:h_0^{\sigma},\ldots,h_n^{\sigma}].
\end{equation}
\end{thm}

\begin{pf*}{Proof.}
Using (\ref{eq:LovExt}), we simply have
\begin{eqnarray*}
\mathbf{E}[g^{(n)}(Y_h)] &=& \int_{[0,1]^n}g^{(n)}[h(\mathbf{x})]\,{\rm d}\mathbf{x}\\
&=& \sum_{\sigma\in\frak{S}_n}\int_{ S_{\sigma}}g^{(n)}\Big[\sum_{i=1}^n \big(h_i^{\sigma}-h_{i-1}^{\sigma}\big)x_{\sigma(i)}\Big]\,{\rm
d}\mathbf{x}.
\end{eqnarray*}
Finally, after an elementary change of variables, we conclude by the Hermite-Genocchi formula (\ref{eq:HermGeno}).\qed
\end{pf*}

Theorem~\ref{thm:ExpgY} provides the expectation $\mathbf{E}[g^{(n)}(Y_h)]$ in terms of the divided differences of $g$ with arguments
$h_0^{\sigma},\ldots,h_n^{\sigma}$ $(\sigma\in\frak{S}_n)$. An explicit formula can be obtained by (\ref{eq:DvDDistArg}) whenever the arguments
are distinct for every $\sigma\in\frak{S}_n$.

Clearly, the special cases
\begin{equation}\label{eq:RmCm}
g(x)=\frac{r!}{(n+r)!}\, x^{n+r},~\frac{r!}{(n+r)!}\, [x-\mathbf{E}(Y_h)]^{n+r},~\mbox{and}~\frac{e^{tx}}{t^n}
\end{equation}
give, respectively, the raw moments, the central moments, and the moment-generating function of $Y_h$. As far as the raw moments are concerned,
we have the following result.

\begin{prop}\label{prop:RM}
For any integer $r\geqslant 1$, setting $A_0:=[n]$, we have,
$$%\begin{equation}\label{eq:rthMoment}
\mathbf{E}[Y_h^r]=\frac{1}{{n+r\choose r}}\sum_{{A_1\subseteq [n]\atop A_2\subseteq A_1}\atop {\cdots\atop A_r\subseteq A_{r-1}}}\,\prod_{i=1}^r
\frac{1}{{|A_{i-1}|\choose |A_i|}}\, v_h(A_i).
$$%\end{equation}
%For any integer $r\geqslant 1$, we have, setting $A_0:=[n]$,
%$$
%\mathbf{E}[Y_h^r]=\frac{1}{{n+r\choose r}}\sum_{A_1\subseteq A_0}\sum_{A_2\subseteq A_1}\cdots\sum_{A_r\subseteq A_{r-1}}\frac{1}{{n\choose
%|A_0\setminus A_1|,|A_1\setminus A_2|,\ldots,|A_r|}}\,\prod_{i=1}^r h(\mathbf{1}_{A_i}).
%$$
%ou encore (mais cette fois $A_0\neq [n]$)
%For any integer $r\geqslant 1$, we have
%$$
%\mathbf{E}[Y_h^r]=\frac{1}{{n+r\choose r}}\sum_{\{A_0,\ldots,A_r\}}\frac{1}{{n\choose
%|A_1|,\ldots,|A_r|}}\,\prod_{i=1}^{r} h(\mathbf{1}_{A_1\cup\cdots\cup A_i})
%$$
%where the sum, having $(r+1)^n$ terms, is taken over all partitions $\{A_0,\ldots,A_r\}$ of $[n]$, with possibly
%empty sets.
\end{prop}

\begin{pf*}{Proof.}
Let $r\geqslant 1$. It can be shown \cite{Ali73} that
$$
\Delta[(\Cdot)^{n+r}:a_0,\ldots,a_n]=\sum_{\textstyle{r_0,\ldots,r_n\geqslant 0\atop r_0+\cdots +r_n=r}} a_0^{r_0}\cdots
a_n^{r_n}=\sum_{0\leqslant i_1\leqslant\cdots\leqslant i_r\leqslant n} a_{i_1}\cdots a_{i_r}.
$$
Hence, from (\ref{eq:ExpgY}) and (\ref{eq:RmCm}) it follows that
\begin{eqnarray*}
\mathbf{E}[Y_h^r]%
&=& \frac{r!}{(n+r)!}\,\sum_{0\leqslant i_1\leqslant\cdots\leqslant i_r\leqslant n}\, \sum_{\sigma\in\frak{S}_n}
h_{i_1}^{\sigma}\cdots h_{i_r}^{\sigma}\\
&=& \frac{r!}{(n+r)!}\,\sum_{0\leqslant i_1\leqslant\cdots\leqslant i_r\leqslant n}\, \sum_{m\in\mathcal{M}_n} v_h(m_{i_1})\cdots v_h(m_{i_r}),
\end{eqnarray*}
where $\mathcal{M}_n$ is the set of the $n!$ maximal chains of the lattice $(2^{[n]},\subseteq)$, and where, for any $m \in \mathcal{M}_n$,
$m_i$ is the unique element of $m$ of cardinality $i$.

For any $B_1 \subseteq \dots \subseteq B_r \subseteq [n]$, let $\mathcal{M}_n^{B_1,\dots,B_r}$ denote the subset of maximal chains of
$(2^{[n]},\subseteq)$ containing $B_1,\dots,B_r$. It is then easy to see that, for any fixed $0 \leqslant i_1 \leqslant \cdots \leqslant i_r
\leqslant n$, the following identity holds:
$$
\bigcup_{\textstyle{B_1 \subseteq \dots \subseteq B_r \subseteq [n] \atop |B_1|=i_1, \dots, |B_r|=i_r}} \mathcal{M}_n^{B_1,\dots,B_r} =
\mathcal{M}_n
$$
and the union is disjoint. Therefore, we have
\begin{eqnarray*}
\mathbf{E}[Y_h^r]&=& %
\frac{r!}{(n+r)!}\,\sum_{0\leqslant i_1\leqslant\cdots\leqslant i_r\leqslant n}\,\sum_{\textstyle{B_1 \subseteq \dots \subseteq B_r \subseteq
[n] \atop |B_1|=i_1, \dots, |B_r|=i_r}}\,\sum_{m\in \mathcal{M}_n^{B_1,\dots,B_r}}
v_h(B_1)\cdots v_h(B_r)\\
&=& \frac{r!}{(n+r)!}\,\sum_{B_1 \subseteq \dots \subseteq B_r \subseteq [n]}|\mathcal{M}_n^{B_1,\dots,B_r}|\, \prod_{i=1}^r v_h(B_i),
\end{eqnarray*}
where
$$
|\mathcal{M}_n^{B_1,\dots,B_r}| = |B_1|!\, (|B_2|-|B_1|)!\, (|B_3|-|B_2|)!\,\cdots (n-|B_r|)!.
$$
Finally, we get the result by setting $A_i:=B_{r+1-i}$ for all $i=1,\ldots,r$.\qed
\end{pf*}

Proposition~\ref{prop:RM} provides an explicit expression for the $r$th raw moment of $Y_h$ as a sum of $(r+1)^n$ terms. For instance, the first
two moments are
\begin{eqnarray*}
\mathbf{E}[Y_h] &=&\frac{1}{n+1}\,\sum_{A\subseteq [n]}\frac{1}{{n\choose |A|}}\, v_h(A),\\
\mathbf{E}[Y_h^2] &=& \frac{2}{(n+1)(n+2)}\,\sum_{A_1\subseteq [n]}\frac{1}{{n\choose |A_1|}}\, v_h(A_1)\sum_{A_2\subseteq
A_1}\frac{1}{{|A_1|\choose |A_2|}}\, v_h(A_2).
\end{eqnarray*}

We now yield a formula for the distribution function $F_{h}(y):=\Pr[Y_h\leqslant y]$ of~$Y_h$.

\begin{thm}\label{thm:CDFYh}
There holds
\begin{equation}\label{eq:CDFYh}
F_{h}(y)=1-\frac {1}{n!}\,\sum_{\sigma\in\frak{S}_n}\Delta[(\Cdot -y)^n_+:h_0^{\sigma},\ldots,h_n^{\sigma}].
\end{equation}
\end{thm}

\begin{pf*}{Proof.}
We have
$$
F_{h}(y) = 1-\Pr[h(\mathbf{X})> y]%
%1-\int_{\{\mathbf{x}\in [0,1]^n\mid h(\mathbf{x})> y\}}\mathrm{d}\mathbf{x}
= 1-\mathbf{E}\big[(Y_h-y)^0_+\big].
$$
Then, using (\ref{eq:ExpgY}) with
$$
g(x)=\frac{1}{n!}\, (x-y)^n_+
$$
leads to the result.\qed
\end{pf*}

It follows from (\ref{eq:CDFYh}) that the distribution function of $Y_h$ is absolutely continuous and hence its probability density function is
simply given by
\begin{equation}\label{eq:fDens}
f_{h}(y)=\frac {1}{(n-1)!}\,\sum_{\sigma\in\frak{S}_n}\Delta[(\Cdot -y)^{n-1}_+:h_0^{\sigma},\ldots,h_n^{\sigma}]
\end{equation}
or, using the B-spline notation (\ref{eq:BSpline}),
$$
f_{h}(y)=\frac {1}{n!}\, \sum_{\sigma\in\frak{S}_n} M(y\mid h_0^{\sigma},\ldots,h_n^{\sigma}).
$$

\begin{rem}\label{rem:2}
\begin{enumerate}
\item[(i)] It is easy to see that (\ref{eq:CDFYh}) can be rewritten by means of the minus truncated power function as
$$
F_{h}(y)=\frac {1}{n!}\,\sum_{\sigma\in\frak{S}_n}\Delta[(\Cdot -y)^n_-:h_0^{\sigma},\ldots,h_n^{\sigma}].
$$
%Particularly, the plus and minus truncated power functions can be used interchangeably in (\ref{eq:fDens});;;;;;;; False! ;;;;;;;;;.

\item[(ii)] When the arguments $h_0^{\sigma},\ldots,h_n^{\sigma}$ are distinct for every $\sigma\in\frak{S}_n$, then combining
(\ref{eq:DvDDistArg}) with (\ref{eq:CDFYh}) immediately yields the following explicit expression
$$
F_{h}(y)=1-\frac{1}{n!}\,\sum_{\sigma\in\frak{S}_n}\sum_{i=0}^n\frac{(h_i^{\sigma}-y)^n_+}{\prod_{j\neq i}(h_i^{\sigma}-h_j^{\sigma})},
$$
or, using the minus truncated power function,
$$
F_{h}(y)=\frac{1}{n!}\,\sum_{\sigma\in\frak{S}_n}\sum_{i=0}^n\frac{(h_i^{\sigma}-y)^n_-}{\prod_{j\neq i}(h_i^{\sigma}-h_j^{\sigma})}.
$$

\item[(iii)] The knowledge of $f_h(y)$ immediately gives an alternative proof of (\ref{eq:ExpgY}). Indeed, using Peano's representation
(\ref{eq:Peano}), we simply have
\begin{eqnarray*}
\mathbf{E}[g^{(n)}(Y_h)] &=& \int_{\R}g^{(n)}(y)\, f_h(y)\,\mathrm{d}y\\
&=& \frac 1{n!}\sum_{\sigma\in\frak{S}_n}\int_{\R} g^{(n)}(y)\, M(y\mid h_0^{\sigma},\ldots,h_n^{\sigma})\,\mathrm{d}y\\
&=& \sum_{\sigma\in\frak{S}_n}\Delta[g:h_0^{\sigma},\ldots,h_n^{\sigma}].
\end{eqnarray*}

\item[(iv)] The case of linear combinations of order statistics is of particular interest. In this case, each $h_i^{\sigma}$ is independent of
$\sigma$ (see Remark~\ref{rem:1}), so that we can write $h_i:=h_i^{\sigma}$. The main formulas then reduce to (see for instance \cite{AdeSan06}
and \cite{AgaDalSin02})
\begin{eqnarray*}
\mathbf{E}[g^{(n)}(Y_h)] &=& n!\, \Delta[g:h_0,\ldots,h_n],\\
F_{h}(y) &=& \Delta[(\Cdot -y)^n_-:h_0,\ldots,h_n],\\
f_{h}(y) &=& M(y\mid h_0,\ldots,h_n).
\end{eqnarray*}
We also note that the Hermite-Genocchi formula (\ref{eq:HermGeno}) provides nice geometric interpretations of $F_h(y)$ and $f_{h}(y)$ in terms
of volumes of slices and sections of canonical simplices (see also \cite{Ali73} and \cite{Ger81}).
\end{enumerate}
\end{rem}

Both functions $F_h(y)$ and $f_h(y)$ require the computation of divided differences of truncated power functions. On this issue, we recall a
recurrence equation, due to de Boor~\cite{deB72} and rediscovered independently by Varsi~\cite{Var73} (see also \cite{Ali73}), which allows to
compute $\Delta[(\Cdot -y)^{n-1}_+:a_0,\ldots,a_n]$ in $O(n^2)$ time.

Rename as $b_1,\dots,b_r$ the elements $a_i$ such that $a_i < y$ and as $c_1,\dots,c_s$ the elements $a_i$ such that $a_i \geqslant y$ so that
$r+s=n+1$. Then, the unique solution of the recurrence equation
\begin{equation}\label{eq:deBoor}
\alpha_{k,l} = \frac{(c_l - y) \alpha_{k-1,l} + (y - b_k) \alpha_{k,l-1}}{c_l - b_k} \qquad (k\leqslant r,\, l\leqslant s),
\end{equation}
with initial values $\alpha_{1,1} = (c_1 - b_1)^{-1}$ and $\alpha_{0,l} = \alpha_{k,0} = 0$ for all $l,k \geqslant 2$, is given by
$$
\alpha_{k,l} := \Delta[ (\Cdot-y)_+^{k+l-2} : b_1,\dots,b_k,c_1,\dots,c_l]\qquad (k + l \geqslant 2).
$$
In order to compute $\Delta[(\Cdot -y)^{n-1}_+:a_0,\ldots,a_n] = \alpha_{r,s}$, it suffices therefore to compute the sequence $\alpha_{k,l}$ for
$k+l \geqslant 2$, $k \leqslant r$, $l \leqslant s$, by means of two nested loops, one on $k$, the other on $l$.

We can compute $\Delta[(\Cdot -y)^n_-:a_0,\ldots,a_n]$ similarly. Indeed, the same recurrence equation applied to the initial values
$\alpha_{0,l} = 0$ for all $l \geqslant 1$ and $\alpha_{k,0} = 1$ for all $k \geqslant 1$, produces the solution
$$
\alpha_{k,l} := \Delta[ (\Cdot-y)_-^{k+l-1} : b_1,\dots,b_k,c_1,\dots,c_l]\qquad (k + l \geqslant 1).
$$
See for instance \cite{Ali73} and \cite{Var73} for further details.

\section{Application to aggregation theory}

As we have already mentioned, the concept of linear combination of lattice polynomials, when it is nondecreasing in each variable, is known in
aggregation theory as the discrete {\em Choquet integral}, which is extensively used in non-additive expected utility theory, cooperative game
theory, complexity analysis, measure theory, etc.\ (see \cite{GraMurSug00} for an overview.) For instance, when a discrete Choquet integral is
used as an aggregation tool in a given decision making problem, it is then very informative for the decision maker to know its distribution. In
that context, the most natural {\em a priori}\/ density on $[0,1]^n$ is the uniform one, which makes the results derived here of particular
interest.

\begin{exmp}
Let $h:[0,1]^3\to\R$ be a discrete Choquet integral defined by $v_h(\{1\}) = 0.1$, $v_h(\{2\}) = 0.6$, $v_h(\{3\}) = v_h(\{1,2\}) = v_h(\{1,3\})
= v_h(\{2,3\}) = 0.9$, and $v_h(\{1,2,3\}) = 1$. According to (\ref{eq:MoEx}), it can be written as
\begin{eqnarray*}
h(\mathbf{x})&=& 0.1\, x_1+0.6\, x_2 + 0.9\, x_3\\
&& \null + 0.2(x_1\wedge x_2)-0.1(x_1\wedge x_3)-0.6(x_2\wedge x_3)\\
&&\null -0.1(x_1\wedge x_2\wedge x_3).
\end{eqnarray*}
Its density, which can be computed through (\ref{eq:fDens}) and the recurrence equation~(\ref{eq:deBoor}), is represented in
Figure~\ref{fig:density} by the solid line. The dotted line represents the density estimated by the kernel method from 10 000 randomly generated
realizations. The typical value and standard deviation can also be calculated through the raw moments: we have
$$
\mathbf{E}[Y_h]\approx 0.608 \quad\mbox{and}\quad \sqrt{\mathbf{E}[Y_h^2]-\mathbf{E}[Y_h]^2}\approx 0.204.
$$
\begin{figure}[t]
\begin{center}
\includegraphics*[width=0.8\linewidth]{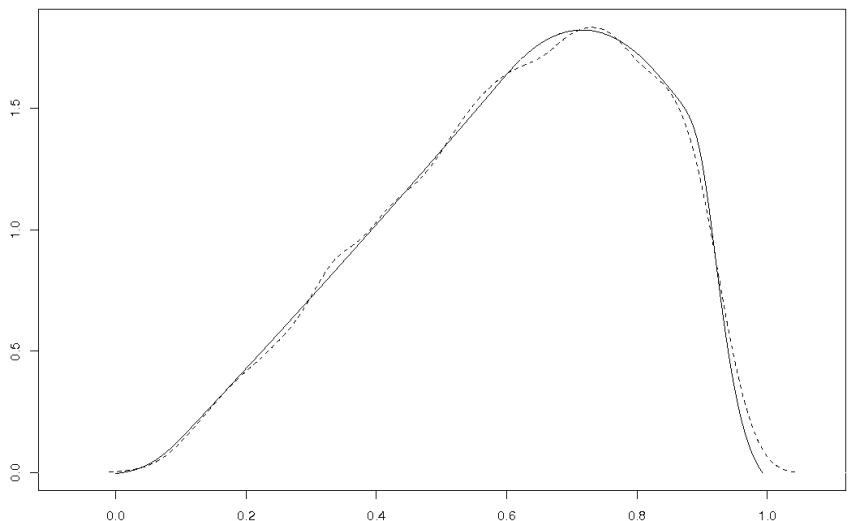}
\caption{\label{fig:density} Density of a discrete Choquet integral (solid line).}
\end{center}
\end{figure}
\end{exmp}

%\bibliographystyle{abbrv}   % styles: plain, unsrt, alpha, abbrv, ieeetr, acm, siam, apalike, amsplain,...
%\bibliography{ReferencesMarichal}

\end{document}